\documentclass[12pt]{article}  
\usepackage{times,amsmath,amsfonts,amsthm,amssymb,eucal}
\usepackage{array}
\usepackage{mathrsfs}
\usepackage{amsmath}
\usepackage{setspace}
\usepackage{fullpage}
\usepackage{fancyhdr}
\usepackage{setspace}
\usepackage{amssymb}
\usepackage{epsfig}
\usepackage{graphicx}
\usepackage{makeidx}
\usepackage{color}
\usepackage{multirow}
\usepackage{epsfig}
\usepackage{graphicx}
\usepackage{makeidx}
\usepackage{gensymb}
\usepackage{enumerate}
\usepackage{multirow}
 \usepackage{setspace}
\usepackage{amssymb}
 \usepackage{setspace}
 \usepackage{mathrsfs}
\usepackage{mathtools}

\newtheorem{theorem}{Theorem}
\newtheorem{proposition}{Proposition}
\newtheorem{lemma}{Lemma}
\newtheorem{corollary}{Corollary}
\newtheorem{definition}{Definition}
\def\Theorem{\begin{theorem}\sl}
\def\EndTheorem{\end{theorem}}
\def\Proposition{\begin{proposition}\sl}
\def\EndProposition{\end{proposition}}
\def\Lemma{\begin{lemma}\sl}
\def\EndLemma{\end{lemma}}
\def\Corollary{\begin{corollary}\sl}
\def\EndCorollary{\end{corollary}}
\def\Definition{\begin{definition}\sl}
\def\EndDefinition{\end{definition}}
\theoremstyle{remark}

\numberwithin{equation}{section}

\begin{document}
\title{ \textbf{On Approximations of the Beta Process in Latent Feature Models}}   

\author{Luai Al Labadi \thanks{{\em Address for correspondence}: Luai Al Labadi, Department of Statistical Sciences, University of Toronto, Toronto, Ontario M5S 3G3, Canada. E-mail: luai.allabadi@utoronto.ca.} \& Mahmoud Zarepour\thanks{M. Zarepour, Department of Mathematics and Statistics,
University of Ottawa, Ottawa, Ontario, K1N 6N5, Canada. E-mail: zarepour@uottawa.ca.} }        

\date{\large{\today}}    
\maketitle
\pagestyle {myheadings} \markboth {} {Approximations of the Beta Process}
\begin{abstract}
The beta  process has recently been widely used as a nonparametric prior for different models in machine learning, including latent feature models. In this paper, we prove the asymptotic consistency of the finite dimensional approximation  of the beta process due to Paisley \& Carin (2009). In addition, we derive an almost sure approximation of the beta process.  This approximation provides a direct method to  efficiently simulate the beta  process.  A simulated example, illustrating the work of the method and comparing its performance to several existing algorithms, is also included.


\vspace{9pt}

\end{abstract}


\section{Introduction}
\label{intro}
The beta process was introduced by Hjort (1990) and later clarified  by Kim (1999) in the field of Bayesian survival analysis as a class of prior processes for cumulative hazard function. Most recently, Kim et al. (2013) used the beta process to simulate the beta-Dirichlet process, a nonparametric prior
for the cumulative intensity functions of a Markov process. Developments of the beta process in machine learning first appeared in the work of Thibaux \& Jordan (2007), where it was shown through the application of document classification that the beta process could be used as a nonparametric prior in latent feature models. They also demonstrated that, when the beta process is marginalized out, one can obtain the Indian buffet process first defined in Griffiths \&  Ghahramani  (2006).  Since then, the beta process has been considered for many other applications in machine learning including factor analysis (Paisley \& Carin, 2009), featural representations of multiple time series (Fox et al., 2009), Gene-expression analysis (Chen et al., 2010),  linear regression (Chen et al., 2010), dictionary learning for image processing (Zhou et al., 2011), and image interpolation (Zhou et al., 2012).

Deriving a  stick breaking  representation for the beta process started by the work of Teh et al. (2007), who developed the stick breaking  representation  of  the Indian buffet process. Recently, a stick-breaking construction of the  full beta process was derived by Paisley et al. (2010). The derivation relied on  a limiting process involving finite matrices, analogous to the limiting process used to derive the Indian buffet process. Broderick et al. (2012) demonstrated that the stick-breaking construction of the beta process can be directly obtained from the characterization of the beta process as a Poisson process. A finite approximation of the beta process  was  suggested without proof by Paisley  \& Carin (2009). To this date, there is no mathematical proof for this approximation despite its use in several applications, including those previously mentioned. Providing a precise proof for the finite approximation of the beta process is the first goal of this paper.

Sampling from the beta process plays a central role in applications including latent feature models. For example, in factor analysis models (West, 2003; Paisley \& Carin, 2009;  Broderick et al., 2012),  the data matrix is decomposed into the product of two matrices plus noise. The model takes the form:
$$X=Z \Phi+ E,$$
where $X \in \mathbb{R}^{N\times P}$ is the data matrix and $E \in \mathbb{R}^{N\times P}$ is an error matrix. The matrix $\Phi\in \mathbb{R}^{K\times P}$ is a matrix of factors, and $Z \in \mathbb{R}^{N\times K}$ is a binary matrix of factor loadings. The dimension $K$ is infinite, and thus the rows of $\Phi$ consist of  an infinite collection of factors. The matrix $Z$ is formed  via a draw from a beta-Bernoulli process. First a sample  from the beta process is drawn. Then  applying this draw to a Bernoulli process yields an infinite binary vector of the matrix $Z$.  The previous step is repeated to generate the matrix $Z$, where each successive draw of the Bernoulli process yields a further row of $Z$.  In other words, the beta process is used to provide an infinite collection of coin-tossing probabilities. Tossing these coins corresponds to a draw  from the Bernoulli process, yielding an infinite binary vector that is considered as a latent feature vector (Broderick et al., 2012). Deriving a  simple, yet efficient, way to simulate the beta  process is the second contribution for this paper.


The rest of the paper is organized as follows. In Section 2, we introduce   the beta process and  its conjugate process, the Bernoulli process. In Section 3, we prove the finite dimensional approximation  of the beta process.  In Section 4, an efficient and convenient method for simulating the beta process is proposed. The approach is based on deriving a finite sum-representation which converges almost surely to the  Ferguson \&  Klass representation (1972) of the beta process. An example illustrating the method and its performance  to other  existing approximations is presented in Section 5.  Finally, our findings are briefly summarized in Section 6.

\section{The Beta Process and the Bernoulli Process}
The beta process and the Bernoulli process are examples of a general family of random measures known as \emph{completely random measures}. Consider a space $\mathbb{X}$ with a $\sigma-$algebra  $\mathscr{B}$  of subsets of $\mathbb{X}$.  A random measure $\Phi$ is said to be completely random measure if
for any finite collection  $A_1, \ldots, A_n$ of disjoint members of $\mathscr{B}$, the random variables $\Phi(A_1),\ldots, \Phi(A_n)$ are independent.
For more details about  completely random measures, consult Kingman (1967). Let $B_0$  be a fixed continuous (non-atomic)
finite measure on $(\mathbb{X},\mathscr{B})$  and $c$ be a positive number. Following Thibaux \& Jordan (2007), the  beta process
$B$, written $B \sim BP \left(c,B_0\right)$, is a completely random measure with L\'evy measure
 \begin{equation}
\nu(d\omega,ds)=cs^{-1}(1-s)^{c-1}dsdB_{0}(\omega), \ \ 0<s<1, \omega \ge 0. \label{eq1}
\end{equation}
For any $S \in \mathbb{X}$, we have (Hjort, 1990; Thibaux \& Jordan, 2007):
 \begin{equation}
 E\left[B(S)\right]=B_0(S)  \quad  \text{and} \quad Var \left[B(S)\right]=\frac{B_0(S)}{c+1}. \label{var}
  \end{equation}
  As in the Dirichlet process (Ferguson, 1973), $c$ is called the \emph{concentration parameter} and  $B_0$ is
  called the \emph{base measure}. Note that, in general,  $c$ can be a positive function of $\omega$, but this is not
  commonly used in latent feature models.  The total mass of $B_0$, $\gamma:=B_0(\mathbb{X})$, is called the \emph{mass parameter}. A draw
  $B \sim \text{BP}\left(c,B_0\right)$ is described by:
 \begin{equation}
 B=\sum_{i=1}^{\infty} p_i\delta_{\omega_i},\label{eq2}
\end{equation}
where $(p_1,\omega_1),(p_2,\omega_2),\ldots$ are the set of atoms in a realization of a nonhomogeneous Poisson process with mean measure $\nu$. Here and throughout the present  paper, $\delta_X$ denotes the Dirac measure at $X$, i.e. $\delta_X(A)=1$ if $X \in A$ and $0$ otherwise for a set $A \in \mathscr{B}$.
As shown in (\ref{eq2}), $B$ is a discrete random measure (with probability 1). Note that, $B$ is a finite measure since
$E[\sum_{i=1}^{\infty}p_i]=E\left[B(\mathbb{X})\right]=B_0(\mathbb{X})=\gamma$.

The stick-breaking representation of the beta process takes the form: (Paisley et al., 2010; Broderick et al., 2012)
\begin{equation}
B=\sum_{i=1}^{\infty} {\sum_{j=1}^{C_i} V^{(i)}_{i,j}\prod_{l=1}^{i-1}(1-V^{(l)}_{i,j})\delta_{\omega_{i,j}}},\label{stickbeta}
\end{equation}
where $C_i \overset{i.i.d.} \sim Poisson(\gamma)$, $V^{(l)}_{i,j}\overset{i.i.d.} \sim beta(1,c)$, and $\omega_{i,j}\overset{i.i.d.} \sim B_0/\gamma$. When $C_i=0$, then the corresponding sum is taken to be zero. The key difference between the stick-breaking representation of the Dirichlet process (Sethuraman, 1994)
and that of the beta process is the weights (probabilities). The weights  in the Dirichlet process depend on each other, while this is not the case for the beta process. Specifically,  the weights that result from  the stick-breaking representation of the Dirichlet process all come from a single stick (the unit interval). Thus, they add up to one. On the other hand, in the beta process, the weights all come from different unit intervals. So they need not add to one. However, as pointed out previously, their sum is  finite almost surely (a.s.).


A connection between the beta process and the Poisson process was established in Paisley et al. (2012) and Broderick et al. (2012). They demonstrated that the beta process  is  a Poisson process with the same mean measure (\ref{eq1}). In particular, Paisley et al. (2012) showed the the stick-breaking construction defined in (\ref{stickbeta}) is equivalent to
\begin{equation}
B\overset{d}=\sum_{j=1}^{C_1}V_{1,j}\delta_{\omega_{1,j}}+\sum_{i=2}^{\infty} {\sum_{j=1}^{C_i} V_{i,j}e^{-T_{ij}}\delta_{\omega_{i,j}}},\label{Poissonbeta}
\end{equation}
where $V_{i,j}\overset{i.i.d.} \sim beta(1,c)$, $T_{ij}\overset{ind.} \sim  gamma(i-1,c)$, $C_i$ and $\omega_{i,j}$ are as defined in (\ref{stickbeta}). In this paper, ``$\overset{d} =$", ``$\overset{d} \to$", ``$\overset{v} \to$" and  ``$\overset{a.s.} \to$"  denote  equal in distribution, convergence in distribution, vague convergence and  almost sure convergence, respectively. More details about convergence of random measures  given in Appendix A. In addition, we use the same notation for the probability measure and  its corresponding cumulative distribution function, i.e. $B(t)=B\left((-\infty,t]\right)$ for $t \in \mathbb{X}=\mathbb{R}$. The inverse of a distribution function (or measure) $B$ is defined by
\begin{equation}
\nonumber B^{-1}(t)=\inf \left\{x:B(x)\ge t \right\}, \ \ 0<t<1.
\end{equation}

A direct link to the Poisson process is the following representation:
\begin{equation}
B\overset{d}=\sum_{i=1}^{\infty} {\sum_{j=1}^{C_i} \left(e^{-\Gamma_{i-1,j}\big/c}-e^{-\Gamma_{i,j}\big/c}\right)\delta_{\omega_{i,j}}},\label{poisson}
\end{equation}
where $\Gamma_{i,j}=E^{(1)}_{i,j}+\cdots+E^{(i)}_{i,j}$ and $(E^{(l)}_{i,j})_{ 1\le l \le i}$ are i.i.d. random variables with the exponential distribution with mean 1, independent of $\omega_{i,j}$. Here we have $\Gamma_{0,j}=0$. The relationship (\ref{poisson}) can be derived from the fact that
\begin{eqnarray*}
e^{-\Gamma_{i-1,j}\big/c}-e^{-\Gamma_{i,j}\big/c}&=&e^{-E^{(1)}_{i,j}\big/c}\ldots e^{-E^{(i-1)}_{i,j}\big/c}\left(1-e^{-E^{(i)}_{i,j}\big/c}\right)\\
&\overset{d}=& (1-V^{(1)}_{i,j})\ldots (1-V^{(i-1)}_{i,j-1}) V^{(i)}_{i,j},
\end{eqnarray*}
since $e^{-E_{1,j}\big/c}\overset{d}=beta(c,1)\overset{d}=1-beta(1,c)$. See Iswaran \& Zarepour (2000)  for an analogous representation of  the Dirichlet process.  Expanding the summation in (\ref{poisson}) for the first values of $i$ gives:
\begin{eqnarray*}
B&\overset{d}=&\sum_{j=1}^{C_1} \left(1-e^{-\Gamma_{1,j}\big/c}\right)\delta_{\omega_{1,j}}\\
&&+\sum_{j=1}^{C_2} \left(e^{-\Gamma_{1,j}/c}-e^{-\Gamma_{2,j}\big/c}\right)\delta_{\omega_{2,j}}\\
&&+\sum_{j=1}^{C_3} \left(e^{-\Gamma_{2,j}/c}-e^{-\Gamma_{3,j}\big/c}\right)\delta_{\omega_{3,j}}+\cdots.
\end{eqnarray*}
Since,   $\Gamma_{i,j}/i \overset{a.s} \to 1$  as $i \to \infty$ (by the strong law of large numbers), the weights becomes negligible for a large value of $i$. This makes representation (\ref{poisson}) useful for simulation purposes through a truncation approach. See Zarepour \&  Al Labadi (2012) and Al Labadi \& Zarepour (2013a) for  further discussion about  truncation procedures.

As pointed out  earlier, the beta process is useful as a parameter for the Bernoulli process. The Bernoulli process can be defined, in general,  for any base measure on $\mathbb{X}$. In our case,  we consider the the base measure to be $B$, where  $B\sim BP(c,B_0)$. Then a Bernoulli process $Y$ with base measure $B,$ written $Y \sim BeP (B),$ is a completely random measure $$Y=\sum_i b_i\delta_{\omega_i},$$
 where $b_i \overset{i.i.d.}\sim Bernoulli (p_i)$ for $p_i$  given in (\ref{eq2}). Observe that, since $E[Y|B]=B$ and $B$ is a finite measure, the number of non-zero points in any realization of the Bernoulli process is finite.

The following theorem shows that the beta process is the conjugate prior for the Bernoulli process. This conjugacy extends the conjugacy between the Bernoulli and beta distributions. For the proof of the theorem, see Thibaux \& Jordan (2007)  and Miller (2011).

\Theorem
\label{BSP3}
Let $B \sim \text{BP}\left(c,B_0\right)$,  and, for $i=1,\ldots,m$, let $X_i|B \sim \text{BeP}(B)$  be $m$ independent Bernoulli process draws from $B.$ The posterior distribution of $B$ after observing $X_1,\ldots, X_m$ is still a beta process:
\begin{equation}
B^{*}_m=B|X_1,\ldots,X_m \sim  \text{BP}\left(c^*, B^{*}_{0,m}\right), \nonumber
\end{equation}with $c^{*}=c+m$ and
\begin{eqnarray}
\nonumber B^{*}_{0,m}&=&\frac{c}{c+m}B_0+\frac{1}{c+m} \sum_{i=1}^{m}X_i.
\end{eqnarray}
\EndTheorem


\section{Finite Dimensional Approximation of the Beta Process}
In this section, we prove  convergence of the finite approximation of the beta process, which  was originally  proposed by  Paisley \&  Carin (2009) without a proof. As mentioned in the Introduction, this approximation plays a crucial role in several applications. 
\Theorem
\label{BSP0}
Consider a space  $\left(\mathbb{R},\mathcal{B}\right)$,  where $\mathbb{R}$ denotes the real line and $\mathcal{B}$ is the Borel $\sigma-$algeba of subsets of $\mathbb{R}$. Let $B_0$ be a finite continuous fixed  measure on  $(\mathbb{R},\mathcal{B})$ with $B_0(\mathbb{R})=\gamma$ and $c$ be a positive number. For $n > \gamma$, define the process $B_n$ as follows:
$$B_n=\sum_{i=1}^np_{i,n}\delta_{\omega_i},$$
$$p_{i,n} \overset{i.i.d.} \sim Beta\left(\frac{c\gamma}{n},c\left(1-\frac{\gamma}{n}\right)\right),$$
$$\omega_{i} \overset{i.i.d.} \sim B_0/\gamma,$$
$$\left(p_{i,n}\right)_{1 \le i \le n} \text { and } \left(\omega_i\right)_{1 \le i \le n} \text { are independent}.$$
Then, as $n\to \infty$, $B_n \overset{d} \to B,$ where $B\sim BP\left(c, B_0\right).$
\EndTheorem

\proof The proof is decomposed in several parts.

\noindent \textbf{Part I:} We apply Proposition \ref{Prop3.21} (Appendix A) to  show that, as $n\to \infty$,
\begin{eqnarray}
 \nonumber nP\left[p_{1,n}\in (x,1)\right]=\frac {n\Gamma(c)}{\Gamma({c\gamma}/{n})\Gamma(c-{c\gamma}/{n})}\int_{x}^1s^{{c\gamma}/{n}-1}(1-s)^{c(1-{\gamma}/{n})-1}ds &&\\
\label{eq6} \overset{v}\to \mu(x)=c\gamma\int_{x}^1 s^{-1}(1-s)^{c-1}ds.
\end{eqnarray}
Observe that, for any $x > 0$, $\Gamma(x)={\Gamma(x+1)}/{x}.$ With $x=c\gamma/n,$ we obtain ${n}/{\Gamma(c\gamma/n)}= {c\gamma}/{\Gamma(c\gamma/n+1)}.$ Since $\Gamma(x)$ is a continuous function, as $n\to \infty$, we get
$$\frac {n}{\Gamma\left({c\gamma}/{n}\right)}\to c\gamma$$
and
$$\frac{\Gamma\left(c\right)}{\Gamma\left(c-{c\gamma}/{n}\right)}\to 1.$$ It follows that, as $n\to \infty$,
$$\frac {n\Gamma(c)}{\Gamma({c\gamma}/{n})\Gamma(c-{c\gamma}/{n})} \to c\gamma.$$
On the other hand,  since $x<s<1$, we have $s^{-1}<x^{-1}$ and $s^{c/n}<1$. Thus, $s^{{c}/{n}-1}<x^{-1}.$ Consequently, the integrand in (\ref{eq6}) is dominated by  $x^{-1}(1-s)^{c(1-{1}/{n})-1},$ which is integrable for $ x<s<1.$ Therefore, by the dominated convergence theorem,   (\ref{eq6}) holds. It follows immediately from  Proposition \ref{Prop3.21}  that, as $n \to \infty$,
\begin{equation}
\xi_n=\sum_{i=1}^{n}\delta_{\left(p_{i,n},\omega_i\right)} \overset{d} \to \xi, \nonumber
\end{equation}
where $\xi$ is a Poisson random measure with  mean $d\mu \times dB_0/\gamma$ and $\omega_{i}$ is as defined in the statement of Theorem 2.

\noindent \textbf{Part II:} We show that
$$\xi=\sum_{i=1}^{\infty}\delta_{\left(\mu^{-1}(\Gamma_i),\omega_i\right)},$$
where  $\Gamma_i=E_1+\cdots+E_i$ and $(E_i)_{i \ge 1}$ is  a sequence of  i.i.d. random variables with the exponential distribution with mean 1.
It is clear that $\sum_{i=1}^{\infty}\delta_{\Gamma_i}$
is a a PRM($\lambda$), where $\lambda$ is Lebesgue measure (Resnick, 2006, Example 5.1). By the special case of Proposition 3.9 of Resnick (1987), $\sum_{i=1}^{\infty}\delta_{\left(\Gamma_i,\omega_i\right)}$ is a PRM($\lambda \times B_0/\gamma$).  Define $T:[0,\infty) \times (-\infty,\infty) \to [0,\infty) \times (-\infty,\infty)$ via $T(x,y)=\left(\mu^{-1}(x),y\right)$. If $t>0$ and $a<b$, we have
\begin{eqnarray*}
 \left(\lambda \times B_0/\gamma\right)\circ T^{-1}\left([t,\infty)\times(a,b)\right)&=&\left(\lambda \times B_0/\gamma\right)\left\{(x,y): \mu^{-1}(x)\ge t \text{ and } a<y<b  \right\}\\
              &=&\left(\lambda \times B_0/\gamma\right)\left\{(x,y): x \le \mu(t) \text{ and } a<y<b \right\}\\
              &=&\lambda\left([0,\mu(t)]\right) B_0 ((a,b))/\gamma\\
              &=&\mu (t)B_0 ((a,b))/\gamma.
\end{eqnarray*}

\noindent Thus, by Proposition 3.7 of Resnick (1987), $\xi=\sum_{i=1}^{\infty}\delta_{\left(\mu^{-1}(\Gamma_i),\omega_i\right)}$ is a PRM($d\mu \times dB_0/\gamma$). Therefore, as $n \to \infty$,
\begin{equation}
\xi_{n}=\sum_{i=1}^{n}\delta_{\left(p_{i,n},\omega_i\right)} \overset{d} \to \xi=\sum_{i=1}^{\infty}\delta_{\left(\mu^{-1}(\Gamma_i),\omega_i\right)}. \label{convergence}
\end{equation}

\noindent \textbf{Part III:} For $h>0$, the map
 $$T_h\left(\sum_{i=1}^{n} \delta_{\left(p_{i,n},\omega_i\right)}\right)=\sum_{i=1}^{n} {{p_{i,n}}I\left({p_{i,n}}>h\right)}\delta_{\omega_i}$$
defined on the set of point processes is continuous with respect to vague topology for random measures (there are finite number of terms in the summation). Therefore, as $n \to \infty$ and for $h>0$, applying $T_h$ to  (\ref{convergence}), we obtain by the continuous mapping theorem (Resnick, 1987, p. 152)
$$B_{n,h}=\sum_{i=1}^{n} {{p_{i,n}}I\left({p_{i,n}}>h\right)\delta_{\omega_i}}\overset{d}\to B_h=\sum_{i=1}^{\infty} {{\mu^{-1}(\Gamma_i)}I\left({\mu^{-1}(\Gamma_i)}>h\right)\delta_{\omega_i}}.$$

\noindent Note that, as $h \to 0$, $B_h \overset{d}\to B$, where $B=\sum_{i=1}^{\infty} {{\mu^{-1}(\Gamma_i)}\delta_{\omega_i}}$. To complete the proof, by Theorem 3.2 of Billingsley (1999), it remains  to show that for any Borel set $A$,
\begin{eqnarray}
\nonumber D_n^{h,\epsilon}(A)&=&P \left\{\left|B_n(A)-B_{n,h}(A)\right| \ge  \epsilon\right\} \\
 \nonumber                         &=&P \left\{\left|\sum_{i=1}^{n} p_{i,n}\delta_{\omega_i}(A)-\sum_{i=1}^{n} {p_{i,n}I\left(p_{i,n}\ge h\right)\delta_{\omega_i}}(A)\right| \ge  \epsilon\right\}\to 0,
\end{eqnarray}

as $n \to \infty$ and $h \to 0$. We have
\begin{eqnarray}
\nonumber \lefteqn{ D_n^{h,\epsilon}(A)=P \left\{\sum_{i=1}^{n} {{p_{i,n}}I\left({p_{i,n}}\le h\right)\delta_{\omega_i}(A)}\ge \epsilon\right\}}\\
\nonumber  & &\le \epsilon^{-1}E\left[\sum_{i=1}^{n}{p_{i,n}} I\left( p_{i,n} \le h\right)\delta_{\omega_i}(A)\right]\\
\nonumber  & &\le \epsilon^{-1}\sum_{i=1}^{n}E\left[{p_{i,n}} I\left( p_{i,n}  \le h\right)\right]E\left[\delta_{\omega_i}(A)\right]\\
\nonumber  & & = \epsilon^{-1}nE\left[p_{1,n}I(p_{1,n}\le h)\right]E\left[\delta_{\omega_1}(A)\right]
\end{eqnarray}

\begin{eqnarray}
 \nonumber  & & =c\epsilon^{-1}E\left[\delta_{\omega_1}(A)\right]\int_0^h x nP\left\{p_{1,n}\in dx\right\}.
\end{eqnarray}

\noindent By (\ref{eq6}), as $n \to \infty,$  we get
\begin{eqnarray}
\label{eq:54}  D_n^{h,\epsilon}(A)  \to c\epsilon^{-1} E\left[\delta_{\omega_1}(A)\right]\int_0^h (1-x)^{c-1}dx.
\end{eqnarray}
Observe that, the integral in (\ref{eq:54})  goes to zero as $h \downarrow 0$. Therefore, by Theorem 3.2 of Billingsley (1999), $B_n \overset{d}\to B,$
where
\begin{equation}
 B=\sum_i \mu^{-1}(\Gamma_i)\delta_{\omega_i}\label{eq20}
\end{equation}
is the  Ferguson \&  Klass  (1972) representation of the beta process. This completes the proof of the theorem.
\endproof


\section{A New Algorithm to Generate the Beta Process}
There are two general techniques to write a series representation for any L\'evy process having no Gaussian component.
The first one is due  to Ferguson \&  Klass (1972). The second technique is due to Wolpert \& Ickstadt (1998). Since we consider only the  case when the concentration parameter is constant (i.e., $c(\omega)=c$ for all $\omega \ge 0$), the approach of Wolpert \&  Ickstadt (1998) and the approach of Ferguson  \&  Klass  (1972) are equivalent. See also Al Labadi \&  Zarepour (2013b, 2014a,b) for further discussion about the two approaches. Therefore, the algorithm of Wolpert \&  Ickstadt (1998) is not included here.

Let $B \sim {BP}\left(c,B_0\right)$ with a continuous $B_0$. Then the Ferguson \&  Klass  (1972) representation of $B$ takes the form given in (\ref{eq20}). Since no closed form for the inverse of the L\'evy measure (\ref{eq1}) exists, working with (\ref{eq20}) is relatively difficult in practice.  The next theorem outlines a remedy to this problem, where an almost sure approximation to (\ref{eq20}) is developed based on a similar result in Zarepour \&  Al Labadi (2012) for the Dirichlet process. Convergence of  random measures is taken with
respect to the vague topology on the space of point measures. Consult Appendix A for a background on convergence of random measures.
\Theorem \label{Beta}
Let $B_0$  be a fixed continuous finite measure on $\mathbb{R}$ with $B_0(\mathbb{R})=\gamma$.  Let $(\omega_i)_{i \geq 1}$ be i.i.d. random variables with common distribution $B_0/\gamma$ and $\Gamma_i=E_1+\cdots+E_i,$ where $\left(E_i\right)_{i\ge1}$ are i.i.d. random variables with the exponential distribution with mean 1, independent of $(\omega_i)_{i\ge1}.$  For $n> \gamma$, define $$\mu_{n}(x)=\frac{\Gamma\left(c\right)}{\Gamma({\gamma c}/{n})\Gamma\left(c-{c\gamma}/{n}\right)}\int_{x}^1
s^{{c\gamma}/{n}-1}\left(1-s\right)^{c\left(1-{\gamma}/{n}\right)-1}ds.$$ Then, as $n\to \infty$,
\begin{equation}
B_{n}=\sum_{i=1}^{n} {{\mu^{-1}_{n}\left(\frac{\Gamma_i}{\Gamma_{n+1}}\right)}\delta_{\omega_i}} \overset{a.s.}\to
B=\sum_{i=1}^{\infty} {\mu^{-1}\left(\Gamma_i\right)\delta_{\omega_i}}. \nonumber
\end{equation}
\EndTheorem

\proof It follows from the proof of Theorem 2 that, for any $x>0$,
\begin{equation}
n \mu_n(x) \to \mu(x). \label{Th4-1}
\end{equation}
Notice that, the left hand side of (\ref{Th4-1}) is a sequence of
a continuous monotone functions converging to a monotone function for every $x>0$. This is
equivalent to the convergence of their inverse function to the inverse
function of the right hand side (Resnick, 1987, Proposition 0.1). Thus,
\begin{equation}
 \mu_n^{-1}\left(\frac{x}{n}\right) \to \mu^{-1}\left(x\right). \label{Th4-2}
\end{equation}
Now, taking $x=\Gamma_i$ in (\ref{Th4-2}) and  the fact that $\Gamma_{n+1}/n \to 1$  as $n\to \infty$
 (by the strong law of large numbers) we get
\begin{equation}
 \mu_n^{-1}\left(\frac{\Gamma_i}{\Gamma_{n+1}}\right) \to \mu^{-1}\left(\Gamma_i\right). \label{Th4-3}
\end{equation}
To prove Theorem 2, by Lemma 1 of Al Labadi and Zarepour (2013b),  we show that, for all $k$ fixed,
\begin{equation}
\sum_{i=1}^{k} {{\mu^{-1}_{n}\left(\frac{\Gamma_i}{\Gamma_{n+1}}\right)}\delta_{\omega_{(i)}}} \overset{a.s.}\to
\sum_{i=1}^{k} {\mu^{-1}\left(\Gamma_i\right)\delta_{\omega_{(i)}}} \nonumber
\end{equation}
as $n \to \infty,$ where   $\omega_{(1)}\le \cdots \le \omega_{(n)}$ represent the corresponding order statistics of $\omega_1,\ldots, \omega_n$.
This directly follows by (\ref{Th4-3}).
\endproof

The next algorithm  is used to generate   samples from an approximation of the  beta process with parameters $c$ and $B_0$, where $B_0$ continuous. The steps of the algorithm are:
\begin{enumerate} [(1)]
\item Fix a relatively large positive integer $n$.
\item Generate  $\omega_i\overset{\text{i.i.d.}}\sim B_0/\gamma$ for $i=1,\ldots,n.$
\item Generate $E_i\overset{\text{i.i.d.}}\sim exponential(1)$  for  $i=1,\ldots,n+1$  such that $\left(E_i\right)_{1\le i \le n+1}$ $\left(\omega_i\right)_{1\le i \le n}$ are independent where $\Gamma_i=E_1+\cdots+E_i.$
\item For $i=1, \ldots,n,$ compute $\left(\mu^{-1}_{n}\left({\Gamma_i}/{\left({\Gamma_{n+1}}\right)}\right)\right),$ which is  simply the quantile function of the ${beta}\left(c\gamma/n,c(1-\gamma/n)\right)$ distribution evaluated at  $1-{\Gamma_i}/ \Gamma_{n+1}$.
\item Set $B_{n}$ as  in Theorem \ref{Beta}.
 \end{enumerate}

\section{Empirical Results: Comparison with Other Methods}

Several algorithms  to sample from the beta process   exist in the literature. In this section, we compare the new approximation of the beta process with  the finite dimensional approximation of the beta process (Paisley \& Carin,2009), Damien et al.'s (1995) algorithm, Lee \& Kim's  (2004) algorithm, and Lee's (2007) algorithm. Summary of these algorithms are given in Appendix B. We would like to highlight here that the new approximation converges almost surely to the Ferguson \&  Klass  (1972) representation of the beta process, while the other four algorithms converge only in  distribution. Almost sure convergence is not only stronger than convergence in distribution but also provides a proper path by path comparison between the suggested approximation and the limit. Thus, almost sure convergence is preferred over convergence in distribution.

In order to make comparisons among the  algorithms, we  use equivalent settings for the parameters characterizing these algorithms (see Table 1). We refer the reader to the original papers for the details of the algorithms.  We consider the beta process with $c=2$ and   $B_0(x) = x$ (i.e., a uniform distribution on $[0, 1]$).  We  compute  the absolute maximum difference  between an approximate sample mean and the exact mean. See also Lee \& Kim (2004) and Lee (2007) for similar comparisons.  The exact mean is $x$; see (\ref{var}). We refer to this statistic by the maximum mean error. Specifically,
\begin{eqnarray*}
\text{maximum mean error}=\max_{x}\left|E\left[B_n(x)\right]-E\left[B(x)\right]\right|&=& \max_{x}\left|E\left[B_n(x)\right]-x\right|,
\end{eqnarray*}
where $x=0.1,0.2,\ldots,0.9,1.0$, $B_n$ is an approximation of the beta process, and $B\sim {BP}(c=2,B_0(x)=x).$
Note that $E\left[B_n(x)\right]$ is approximated  by obtaining  the mean  at $x=0.1,0.2,\ldots,0.9,1.0$ of $3000$ i.i.d. sample paths from the approximated process $B_n$. Similarly, we compute the maximum standard deviation  error between an approximate sample standard deviation (s.d.) and the exact standard deviation. The exact standard deviation is $\sqrt{x/3}$; see (\ref{var}). Thus,
\begin{eqnarray*}
\text{maximum s.d. error}=\max_{x}\left|s.d\left[B_n(x)\right]-s.d.\left[B(x)\right]\right|&=& \max_{x}\left|s.d.\left[B_n(x)\right]-\sqrt{x/3}\right|.
\end{eqnarray*}
Table 1 depicts values of the maximum mean  error, the maximum standard deviation  error, and the corresponding computational time. The computational time is computed by applying the code ``System.Time" available in R. As seen in Table 1, the new algorithm has the smallest mean and standard deviation errors. In addition, it  has   a reasonable computation time.

\begin{table}[htbp]
\caption{This table reports the maximum mean error, the maximum standard deviation error, and the corresponding computation time. Here, DSL, KL, and PC stand
for  Damien et al.'s algorithm (1995),  Lee \& Kim's algorithm (2004), and Paisley \&  Carin's (2009) algorithm, respectively.}
\begin{center}
\begin{tabular}{lllll}
\hline
\hline
\multicolumn{0}{c} {Algorithm} &\multicolumn{0}{c} {Parameters} & \multicolumn{0}{c} {max. mean error}  &
\multicolumn{0}{c} {max. s.d. error}& \multicolumn{0}{c} {Time}\\
\hline
DSL &$m=n=200$ & 0.0110& 0.0156&100.70\\
\\
LK  &$\epsilon=0.01$   &   0.0128&0.0107& 0.51\\
\\
Lee  &$n=200, \epsilon=0.05$ & 0.0162&0.0392&1.73\\
\\
PC &$n=200$   & 0.0113 & 0.1280&0.42\\
\\
New &$n=200$   & 0.0087&0.0061&26.38\\
\\
\hline
\end{tabular}
\end{center}
\label{table4.3}
\end{table}

\section{Conclusions}
In this paper, we have proved the finite dimensional approximation  of the beta process (Paisley \& Carin, 2009). This approximation has been used in several  machine learning models. We also have derived an almost surely approximation of the beta process. This new approximation provides a simple, yet efficient, way to simulate the beta  process.


\appendix

\section {Vague Convergence}

The main objective of this appendix is to give a brief introduction about convergence of random measures.   Let $\left(\mathbb{E},\mathscr{E}\right)$  be a state space as before. Let $C_K^+(\mathbb{E})$  be the set of continuous functions $f:\mathbb{E} \to [0,\infty)$ with compact support. A measure $\mu$ is called Radon if $\mu(K)<\infty$ for any compact
set $K$ in $\mathbb{E}$. Let ${M}_+(\mathbb{E})$  be the space of Radon measures in $\mathbb{E}$. Let $\mathscr{M}_+(\mathbb{E})$  be  the smallest
$\sigma-$algebra of subsets of $M_+(\mathbb{E})$ making the maps $\mu \to \mu(f)=\int f(x) d\mu(x)$ from ${M}_+(\mathbb{E})$ to $\mathbb{R}$ measurable for
all functions $f\in C_K^+(\mathbb{E})$. Note that, $\mathscr{M}_{+}(\mathbb{E})$ is the Borel $\sigma-$algebra generated by the topology of vague convergence. If $\mu_n,\mu \in {M}_+(\mathbb{E}),$
we say that $(\mu_n)_n$ converges vaguely to $\mu$ (and we write $\mu_n \overset{v}\to\mu$) if $\mu_n(f) \overset{v}\to\mu(f)$
for any $f\in C_K^+(\mathbb{E}).$

A \emph{random measure} on $\mathbb{E}$ is any measurable map $\xi$ defined on a probability space $(\Omega,\mathscr{A},P)$ with values in
$\left(M_+(\mathbb{E}),\mathscr{M}_+(\mathbb{E})\right)$. If $\xi_n,\xi$ are random measures on $\mathbb{E},$ we say that $(\xi_n)_n$ converges
in distribution to $\xi$ (and we write $\xi_n \overset{d}\to\xi$) if $\left\{P\circ\xi_n^{-1}\right\}_n$ converges weakly to
$P\circ\xi^{-1}.$ By Theorem 4.2 of Kallenberg (1983), $\xi_n \overset {d} \to \xi$ if and only if $\xi_n(f)\to \xi(f),$ i.e.
$$\int_E f(x)\xi_n(dx) \to \int_E f(x) \xi(dx),~ \forall f \in C_K^+(\mathbb{E}).$$

We say that   $(\xi_n)_n$ converges vaguely almost surely to $\xi$ (and write $\xi_n \overset{a.s.}\to\xi$) if there exists a
set $\widetilde{\Omega}\in \mathscr{A}$ with $P(\widetilde{\Omega})=1$ such that
$\forall \omega \in \widetilde{\Omega},$ $\xi_n(\omega,\cdot) \overset {v} \to \xi(\omega,\cdot),$ i.e. $$\int_E f(x)\xi_n(\omega,dx) \to
\int_E f(x) \xi(\omega,dx),~ \forall f \in C_K^+(\mathbb{E}).$$

The space $M_+(\mathbb{E})$ endowed with the vague topology is a complete  separable metric space (Resnick, 1987, Proposition 3.17).  See also Kallenberg (1983).

The next proposition is fundamental in studying the weak convergence of Poisson processes. It gives necessary and sufficient conditions for empirical measures to converge to a Poisson random measure. The proof follows by mimicing the proof of Proposition 3.21 of Resnick (1987) with $j \ge 1$ and $j/n$ are replaced by $1\le j \le n$ and $\omega_j$, respectively. See also Ishwaran, James, and Zarepour (2009) for an analogous result.

\Proposition \label{Prop3.21} Suppose for each $n \ge 1$, we have $\left(X_{j,n}\right)_{1\le j\le n}$ are i.i.d. random variable random elements of $\left(\mathbb{E},\mathscr{E}\right)$ and $\mu$ is a Radon measure on $\left(\mathbb{E},\mathscr{E}\right)$. Let  $\omega_{i} \overset{i.i.d.} \sim F$  such that $(\omega_{j})_{1\le j\le n}$  and $\left(X_{j,n}\right)_{1\le j\le n}$ are independent. Define $\xi_n=\sum_{j=1}^{n}\delta_{\left(\omega_j,X_{j,n}\right)}$ and suppose $\xi$ is a Poisson random measure  on $\mathbb{R}\times \mathbb{E}$ with mean $dF \times d\mu.$ Then $$\xi_n \overset{d} \to \xi  \text{ in } M_p\left(\mathbb{R}\times \mathbb{E}\right)\quad \text{if and only if} \quad nP\left[X_{1,n}\in \cdot\right]\overset{v} \to \mu \quad \text {on } \mathbb{E}.$$
\EndProposition

\section {Other Sampling Algorithms} Several algorithms are suggested to sample from the  beta process $B \sim {BP}\left(c,B_0\right)$ with a continuous $B_0$. We consider  the algorithm of Damien, Laud, and Smith (1995), the algorithm of Lee and Kim  (2004) and  the algorithm Lee (2007). Below is a brief discussion  of these algorithms. We refer the reader to the original papers for more details.

\vspace{3mm}
\noindent \textbf{ $\bullet$ Damien-Laud-Smith Algorithm:} Using  the fact that the distributions of the increments of a nondecreasing L\'evy process
are infinitely divisible, Damien, Laud, and Smith (1995) derived an algorithm to generate approximations for infinitely divisible random
variables and used it to generate  the beta process. Let $p_i$ denotes the increment of the process $B$ in  the interval $\Delta_i=(\omega_{i-1},\omega_i],$ i.e. $p_i=B(\omega_i)-B(\omega_{i-1}).$ The steps of the Damien-Laud-Smith algorithm for simulating  an approximation for the jump $p_i$ are:
\begin{enumerate}
\item [(1)] Fix a  relatively large positive integer $n$.
\item [(2)] Generate  independent values $z_{ij}$ from the probability density function ${dB_{0}(t)}/$ ${B_{0}(\Delta_i)}$, for $j=1,\ldots,n.$
\item [(3)] Generate $x_{ij} \sim {beta} (1,c),$ for $j=1,\dots,n.$
\item [(4)] Generate $y_{ij}$: $y_{ij}|x_{ij} \sim {Poisson}(\lambda_in^{-1}x_{ij}^{-1})$, for $j=1,\dots,n,$ where $\lambda_i=B_0(\omega_i)-B_0(\omega_{i-1})$.
\item [(5)] Set $p_{i,n}=\sum _{j=1}^n x_{ij}y_{ij}.$ For large $n$, $p_{i,n}$ is an approximation of $p_i$.
\end{enumerate}

Damien, Laud, and Smith (1995) showed that  $p_{i,n} \overset{d} \to p_i$, as $n \to \infty$. That is, $p_{i,n}$ is an approximate  sample from the ith increment of $B$. Note that, the Damien-Laud-Smith  algorithm generates only the increments of the process and not the entire process. To obtain the whole process, we set $$B_{m,n}=\sum_{i=1}^m p_{i,n}\delta_{\omega_i}.$$ For large $m$ and $n$,  $B_{m,n}$ is an approximation of $B.$

\vspace{3mm}

\noindent \textbf{$\bullet$ Lee-Kim Algorithm:} The Kim and Lee algorithm for the beta process with parameters $c$ and $B_{0}$ with  $B_0$ continuous can be described as follows. First the L\'evy measure $\nu$ of the beta process given by
 (\ref{eq1}) is approximated by:
\begin{equation}
\nu_{\epsilon}(d\omega,ds)=\frac{c}{\epsilon} b(s:\epsilon,c)dB_{0}(\omega)ds, \nonumber
\end{equation}
where
\begin{equation}
b(x:a,b)=\frac{\Gamma(a+b)}{\Gamma(a)\Gamma(b)}x^{a-1}(1-x)^{b-1}, \quad \text{for } 0<x<1, a>0, b>0.\label{B7}
\end{equation}

The steps of the Lee-Kim algorithm for the beta process $B$ are:
\begin{enumerate}
\item [(1)] Fix a  relatively small positive number $\epsilon$.
\item [(2)] Generate the total number of jumps  $n\sim {Poisson} \left(c\gamma/\epsilon\right)$.
\item [(3)] Generate i.i.d. random variables $\omega_1,\ldots,\omega_n$ from the probability density function $dB_0/\gamma.$
\item [(4)] Let $\omega_{(1)} \le \ldots \le \omega_{(n)}$ be the corresponding order statistics of $\omega_1, \ldots, \omega_n.$
\item [(5)] Generate the jump sizes $p_1,\ldots,p_n:$ $p_i|\omega_{(i)} \sim {Beta}(\epsilon,c).$
\item [(6)] Set $B_{\epsilon}=\sum_{i=1}^n p_i \delta_{\omega_{(i)}}.$
\end{enumerate}
 \vspace{3mm}
Lee and Kim (2004) showed that $B_{\epsilon} \overset{d}\to B,$ as $\epsilon \to 0$.

 \vspace{3mm}

\noindent \textbf{$\bullet$ Lee  Algorithm:} The steps of the Lee algorithm are:

\begin{enumerate}
\item [(1)] Fix a  relatively large positive integer $n$.
\item [(2)] Generate i.i.d. random variables $\omega_1,\ldots,\omega_n$ from the probability density function $dB_0/\gamma.$
\item [(3)] For $i=1,\dots,n,$ generate $x_i \sim b(s:\epsilon,c),$ where $b(s:\epsilon,c)$ is defined in (\ref{B7})
\item [(4)] For $i=1,\dots,n,$ generate $y_i\sim Poisson\left({\gamma b(x_i:1,c)}/{n x_i b(x_i:\epsilon,c)}\right).$
\item [(5)] Set $B_{n}=\sum_{i=1}^n x_i y_i \delta_{\omega_i}.$
\end{enumerate}
\vspace{3mm}

Lee (2007) proved that, as $n\to \infty,$  $B_{n} \overset{d} \to B$.

\end{document}